\documentclass[reqno]{amsart}
\usepackage{epsfig}

\numberwithin{equation}{section}

\newtheorem{lemma}{Lemma}[section]
\newtheorem{prop}[lemma]{Proposition}

\def\SL{\operatorname{SL}}
\def\SO{\operatorname{SO}}
\def\Spin{\operatorname{Spin}}

\def\G{{\mathrm G}}
\def\tG{{\widetilde{\mathrm G}}}
\def\bG{{\overline{\mathrm G}}}
\def\CG{{\mathrm C}}
\def\bP{{\overline{P}}}

\def\im{{\mathrm{Im}\:\delta(\G(k))}}

\def\H{{\mathrm H}}

\def\Hy{{\mathcal H}}
\def\cO{{\mathcal O}}

\def\Q{{\mathbb Q}}

\def\A{{\mathbb A}}

\begin{document}
\title{Addendum to: On Volumes of Arithmetic Quotients of $\SO(1,n)$.}
\author{Mikhail Belolipetsky}
\address{Department of Mathematical Sciences, Durham University, Durham DH1 3LE, UK
\newline\phantom{tt}
Sobolev Institute of Mathematics, Koptyuga 4, 630090 Novosibirsk, Russia}

\email{Mikhail.Belolipetsky@durham.ac.uk}
\subjclass{11E57 (primary); 22E40 (secondary)}
\date{}

\begin{abstract}
There are errors in the proof of the uniqueness of arithmetic subgroups of the
smallest covolume. In this note we correct the proof, obtain certain results which
were stated as a conjecture, and we give several remarks on further developments.
\end{abstract}

\maketitle

\setcounter{section}{1}

\subsection{}
Let us recall some notations and basic notions. Following \cite{B} we will assume
that {\it $n$ is even} and $n\geq4$. The group of orientation
preserving isometries of the hyperbolic $n$-space is isomorphic to
$\SO(1,n)^o$, the connected component of identity of the special orthogonal
group of signature $(1,n)$, which can be identified with
$\SO_0(1,n)$, the subgroup of $\SO(1,n)$ preserving the upper half space.
This group is not Zariski closed in $\SL_{n+1}$ thus in order to
construct arithmetically defined subgroups of $\SO(1,n)^o$ we consider arithmetic
subgroups of the orthogonal group $\SO(1,n)$ or, more precisely, of groups $\G = \SO(f)$
where $f$ is an admissible quadratic form defined over a totally real number field
$k$ (see \cite[Section~2.1]{B}).

We have an exact sequence of $k$-isogenies:
\begin{equation}\label{eq1}
1 \to \CG \to \tG \stackrel{\phi}\to \G \to 1,
\end{equation}
where $\tG(k)\simeq\Spin(f)$ is the simply connected cover of $\G$ and
$\CG\simeq\mu_2$ is the center of $\tG$. This induces an exact sequence in
Galois cohomology (see \cite[Section~2.2.3]{PR})
\begin{equation}\label{eq2}
\tG(k) \stackrel{\phi}\to \G(k) \stackrel{\delta}\to \H^1(k, {\CG}) \to \H^1(k,
\tG).
\end{equation}
The main idea of this note is that by using $(\ref{eq2})$ certain questions about arithmetic
subgroups of $\G$ can be reduced to questions about the Galois cohomology group $\H^1(k, {\CG})$.

A coherent collection of parahoric subgroups $P = (P_v)_{v\in
V_f}$ of $\tG$ ($V_f = V_f(k)$ denotes the set of finite places of the field
$k$) defines a principal arithmetic subgroup $\Lambda = \tG(k)\cap\prod_{v\in
V_f} P_v \subset \tG(k)$ (see \cite{BP}). We fix an infinite place $v$ of $k$
for which $\G(k_v)\simeq\SO(1,n)$ and denote it by $Id$.  The image of $\Lambda$ under
the central $k$-isogeny $\phi$ is an arithmetic subgroup of $\G$ and every maximal
arithmetic subgroup of $\G(k_{Id})$ can be obtained as a normalizer of some
$\phi(\Lambda)$ \cite[Proposition~1.4]{BP}. We will also consider the local
stabilizers of $P$ in the adjoint group $\G (=\bG)$, defining $\bP_v$ to be the
stabilizer of $P_v$ in $\G(k_v)$ and $\bP = (\bP_v)_{v\in V_f}$. Clearly,
$\bP_v\supset \phi(P_v)$. In the notation of \cite{B} the subgroups $\phi(P_v)$ are
called parahoric subgroups of $\G$, however this terminology is non-standard
and we will avoid using it here.
\medskip

\noindent
{\sc Acknowledgements.} This article was written in Summer 2006 while I have been visiting MPIM in Bonn.
I would like to thank Vincent Emery and Gopal Prasad for helpful discussions.

\subsection{}
Given a totally real number field $k$ with the group of units $U$, let
$$
k_\infty^* = \{ a \in k^* \mid a_v > 0 {\rm\ for\ } v\in V_\infty\setminus Id \ \},\ \
U_\infty = U\cap k_\infty^*.
$$

\begin{lemma} \label{lemma} ${\rm Im}(\delta) \simeq k_\infty^*/(k^*)^2$.
\end{lemma}
\begin{proof} From $(\ref{eq2})$ we have
${\rm Im}(\delta\!:\!\G(k)\to\H^1(k,\mu_2)) = {\rm Ker}(\H^1(k,\mu_2) \to \H^1(k,\tG))$.
The Hasse principle for simply connected $k$-groups implies that $\H^1(k,\tG)$
is isomorphic to $\prod_{v\in V_\infty}\H^1(k_v,\tG)$ \cite[Theorem~6.6]{PR}, and
hence
\begin{equation*}
{\rm Im}(\delta) = {\rm Ker}(\H^1(k,\mu_2) \to {\textstyle\prod_{v\in V_\infty}}\H^1(k_v,\tG)).
\end{equation*}
The group $\H^1(k,\mu_2)$ is canonically isomorphic to $k^*/(k^*)^2$ \cite[Lemma~2.6]{PR}.
It is well known that for all $v\in V_\infty$ such that the group $\G(k_v)$ is anisotropic,
the map $\phi$ in $(\ref{eq2})$ is surjective and hence for all such $v$, ${\rm Im}(\delta_v)
= {\rm Ker}(\H^1(k_v,\mu_2) \to \H^1(k_v,\tG))$
is trivial. For the remaining one infinite place $v (=Id)\in V_\infty$, $\phi(\tG(k_v))$ is
a subgroup of index $2$ in $\G(k_v)$ which consists of the orthogonal transformations with
the trivial spinor norm. Collecting this information together we obtain the required
isomorphism.
\end{proof}

\subsection{} The proof of the uniqueness part in \cite[Theorem~4.1]{B}
contains errors but the result is correct. We will now give another argument
for it. In order to do so we first establish a more general fact and then apply
it to the cases considered in \cite{B}.

Let $P = (P_v)_{v\in V_f}$ and $P' = (P'_v)_{v\in V_f}$ be two coherent
collections of parahoric subgroups of $\tG$ such that for all $v\in V_f$,
$P'_v$ is conjugate to $P_v$ under an element of $\G(k_v)$. For all but
finitely many $v$, $P_v = P_v'$ hence there is an element $g\in\G(\A_f)$
($\A_f$ denotes the ring of finite ad\`eles of $k$) such that $P'$ is the
conjugate of $P$ under $g$. We have $\bP = \prod_{v\in V_f}\bP_v$ is the
stabilizer of $P$ in $\G(\A_f)$. The number of distinct $\G(k)$-conjugacy
classes of coherent collections $P'$ as above is the cardinality $c(\bP)$ of
$\mathcal{C}(\bP) = \G(k)\backslash \G(\A_f)/\bP$, which is called the class
group of $\G$ relative to $\bP$. The class number $c(\bP)$ is known to be
finite (see {\it e.g.} \cite[Proposition~3.9]{BP}). The following result
can be used for obtaining further information about its value.
\medskip

\begin{prop} Let $\G = \SO(f)$, $\tG = \Spin(f)$ for an admissible quadratic
form $f$ defined over $k$ and let $P = (P_v)_{v\in V_f}$ be a coherent
collection of parahoric subgroups of $\tG$. The class number $c(P)$ divides the
order $h_{\infty,2}$ of a restricted $2$-class group of $k$ given by
\begin{equation*}
h_{\infty,2} = \frac{2^{[k:\Q]-1}h_2}{[U:U_\infty]},
\end{equation*}
where $h_2$ is the order of the $2$-class group of $k$.
\end{prop}
\begin{proof} Recall two isomorphisms (see \cite[Proposition~8.8]{PR}, a minor
modification is needed in order to adjust the statement to our setting but
the argument remains the same):
\begin{eqnarray*}
\G(k)\backslash \G(\A_f)/\bP & \simeq & \G(\A_f)/\bP\G(k);\\
\G(\A_f)/\bP\G(k) & \simeq & \delta_{\A_f}(\G(\A_f)) / \delta_{\A_f}(\bP\G(k)),
\end{eqnarray*}
where $\delta_{\A_f}$ is the restriction of the product map $\prod_v \G(k_v)
\to \prod_v \H^1(k_v,\CG)$ to $\G(\A_f)$.

For every finite place $v$, $\H^1(k_v,\tG)$ is trivial (see
\cite[Theorem~6.4]{PR}) which implies $\delta_v: \G(k_v)\to \H^1(k_v,\CG)$ is
surjective. Thus the image of $\delta_{\A_f}(\G(\A_f))$ can be identified with
the restricted direct product $\prod\nolimits'\H^1(k_v, \CG)$ with respect to
the subgroups $\delta_v(\bP_v)$. Also $\delta_{\A_f}(\G(k))$ naturally identifies
with the image of $\delta(\G(k))$ in $\H^1(k,\CG)$ under the embedding
$\psi:\H^1(k,\CG) \to \prod\nolimits'\H^1(k_v, \CG)$. Hence we have an
isomorphism
\begin{equation*}
\delta_{\A_f}(\G(\A_f)) / \delta_{\A_f}(\bP\G(k)) \simeq
{\textstyle\prod'}\:\H^1(k_v, \CG) / \left({\textstyle\prod_v}\delta_v(\bP_v)\cdot\psi(\im)\right).
\end{equation*}

The group $\H^1(k_v,\mu_2)$ is canonically isomorphic to $k_v^*/(k_v^*)^2$, by
Lemma~\ref{lemma} $\im \simeq k_\infty^*/(k^*)^2,$ so we obtain
\begin{equation*}
\frac{\prod\nolimits'\H^1(k_v, \CG)}{\prod_v\delta_v(\bP_v)\cdot\psi(\im)}
\simeq
\frac{\prod\nolimits' k_v^*/(k_v^*)^2}{\delta_P \cdot k_\infty^*/(k^*)^2}
\simeq
\frac{J_f}{\delta_P\cdot J_f^2 k^*} \cdot \frac{k^*/k_\infty^*}{U/U_\infty},
\end{equation*}
where $J_f$ is the ring of finite id\`eles of $k$ and $\delta_P$ denotes $\prod_v\delta_v(\bP_v)$.

Now, $\#(J_f/J_f^2k^*) = h_2$, the group $k^*/k_\infty^*$ splits as a product of local factors
and $\#(k^*/k_\infty^*) = 2^{[k:\Q]-1}$ (see \cite[Chapter~6]{L}).
This implies the proposition.
\end{proof}

In order to give a precise formula for the class number $c(P)$ one has to analyze
the image of $\prod_v\delta_v(\bP_v)$ in $\prod\nolimits'\H^1(k_v, \CG)$.
Still in many practical cases this appears to be unnecessary. Thus in order to prove
the uniqueness of the minimal hyperbolic orbifolds we need to consider
$k = \Q[\sqrt{5}]$ (in the compact case) and $k = \Q$ (for the non-compact orbifolds).
In both cases $h_2 = h = 1$. For $k = \Q[\sqrt{5}]$, $U/U_\infty = \{U_\infty, \frac{1+\sqrt{5}}{2}U_\infty \}$
and $[U:U_\infty] = 2$, which implies $h_{\infty,2} = 1$. For $k = \Q$, clearly,
$h_{\infty,2} = 1$ as well. So in all the cases $c(P) = 1$ which implies that the
corresponding arithmetic subgroups are defined uniquely up to a conjugation by $g\in\SO(1,n)$.
It is clear that we can always chose $g\in\SO_0(1,n)$ and therefore the smallest orbifolds
constructed in \cite{B} are unique up to an (orientation preserving) isometry.

\subsection{} We now turn to Conjecture~4.1 and its analogue for the
non-cocompact orbifolds in \cite[Section~4.4]{B}. Recall that in \cite{B} the
numbers $N(r)$, $N'(r)$ were defined for every $r\ge 2$ and estimated from
above. These numbers are related to the index of the principal arithmetic
subgroups in their normalizers. We now prove

\begin{prop} For every $r\ge2$, $N (r) = N'(r) = 1$.
\end{prop}
\begin{proof}
Let $\Lambda$ be a principal arithmetic subgroup of $\tG$ which corresponds to a
compact or non-compact hyperbolic $n$-orbifold of the minimal volume,
$\Lambda' = \phi(\Lambda)$ and $\Gamma = N_\G(\Lambda')$.

From \cite[Proposition~2.9]{BP}, which in turn follows from the work of J.~Rohlfs,
using the fact that the center of our group $\G$ is trivial, we obtain:
$$
[\Gamma:\Lambda'] = \#(\H^1(k,\mu_2)_\Theta \cap \delta(\G(k))) = \# {\rm
Im}(\delta : \G(k)\to\H^1(k,\mu_2))_\Theta.
$$
We can identify the image of $\delta$ by Lemma~\ref{lemma} and then compute ${\rm Im}(\delta)_\Theta$
using \cite[Section~5.1]{BP}. The cases we are interested in are

\begin{eqnarray*}
k & = & \Q:\ {\rm Im}(\delta)_\Theta = \bigg\{ k^{*2}, (-1)k^{*2} \bigg\};\\
k & = & \Q[\sqrt{5}]:\ {\rm Im}(\delta)_\Theta = \bigg\{ k^{*2}, \frac{1-\sqrt{5}}{2}k^{*2}\bigg\}.
\end{eqnarray*}
In both cases $[\Gamma:\Lambda'] = \#{\rm Im}(\delta) = 2$. Now it is easy to
see that $\Lambda' = \phi(\Lambda) \subset \SO_0(1,n)$. From the other side there
always exists $g\in \SO(1,n)\setminus\SO_0(1,n)$ which normalizes $\phi(\Lambda)$.
For example take $g = {\rm diag}(-1,-1,1,\ldots,1)$. As in all the cases
under consideration the quadratic form associated to $\Lambda$ is diagonal
\cite[Sections~4.3, 4.4]{B}, $g$ stabilizes $\Lambda$ and clearly
$g\in \SO(1,n)\setminus\SO_0(1,n)$. From these facts it follows
that $\Lambda'$ is a maximal arithmetic subgroup in $\SO_0(1,n)$ and thus
$N(r) ({\rm or}\ N'(r)) = 1$.

\end{proof}

This proposition makes {\it precise} the statements of Theorem~4.1 and 4.4 of
\cite{B}. It also implies that Table~2 of {\it loc.~cit.} gives the precise values of the
covolumes of the smallest $n$-dimensional hyperbolic orbifolds in even dimensions up to
$18$.

One other corollary is that cocompact and non-cocompact arithmetic subgroups of
$\SO(1,2r)^o$ of the smallest covolumes can be obtained as the stabilizers of
certain lattices described in \cite[Section~4.3]{B}. We remark that since the
fields of definition of the groups have class number $1$, the lattices in both
cases are free as $\cO_k$-modules.

\subsection{} Correction: on p.~765, l.~9 one should read ``grow super-exponentially'' instead
of ``grow exponentially''. (It follows from \cite{B} that the Euler characteristic is
bounded from below by $\mathrm{const}\cdot\left(\prod_{i=1}^{r}\frac{(2i-1)!}{(2\pi)^{2i}}\right)^{[k:\Q]}$
which for large enough $r$ is $\ge \mathrm{const}\cdot(2r-1)!$)
\bigskip

We conclude this addendum with a few remarks on related results which appeared
after the paper was published.

\subsection{} In \cite[Section~4.5]{B} we observed that for $r > 2$ the minimal covolume
among the arithmetic lattices in $\SO(1,2r)$ is attained on a non-uniform lattice. This
interesting phenomenon was first discovered by A.~Lubotzky for $\SL_2$ over local
fields of positive characteristic. Recently, in \cite{S} A. Salehi Golsefidy proved
that lattices of minimal covolume in classical Chevalley groups over local fields
of characteristic $p>7$ are all non-uniform. This result gives further support to
a {\bf conjecture} that {\it generically (i.e. for groups of high enough rank or fields
of large enough positive characteristic) the minimal covolume is always attained on a
non-uniform lattice}.

\subsection{} In \cite{CM} M. Conder and C. Maclachlan constructed a compact orientable
hyperbolic $4$-manifold which has Euler characteristic $16$. The previously known
smallest example which was used in order to formulate the main result in \cite[Section~5]{B}
had $\chi = 26$. The construction of \cite{CM} agrees with our Theorem~5.5 and it also allows
us to give a more precise formulation of the theorem:
\medskip

\noindent\normalsize {\bf Theorem $5.5'.$}
{\it
If there exists a compact orientable arithmetic hyperbolic
$4$-manifold $M$ with $\chi(M) \leq 16$, then $M$ is defined over $\Q[\sqrt{5}]$ and has
the form  $\Gamma_M\backslash\Hy^4$ with $\Gamma_M$ being a torsion-free subgroup of index
$7200\chi(M)$ of the group $\Gamma_1$ of the smallest arithmetic hyperbolic $4$-orbifold.}
\medskip


\end{document}